\documentclass[12pt]{amsart}
\usepackage{amscd}
\usepackage{amssymb,amsmath}
\usepackage{hyperref}
\usepackage{mathrsfs}
\usepackage{graphicx}
\usepackage{enumitem}
\usepackage{romannum}
\usepackage{mathtools}
\usepackage{array}
\usepackage[utf8]{inputenc}
\usepackage{listings}
\usepackage{lipsum,eso-pic,xcolor}
\usepackage{lineno}
\usepackage{tikz}
\usetikzlibrary{arrows,decorations.pathmorphing,backgrounds,positioning,fit}  
\usetikzlibrary{positioning}
\hypersetup{
colorlinks=true,
linkcolor=blue,}

\usetikzlibrary{trees}
\usetikzlibrary{arrows}

\def\supp{\operatorname{supp}}

\def\reg{\operatorname{reg}}
\def\deg{\operatorname{deg}}

\def\gcd{\operatorname{gcd}}
\def\max{\operatorname{max}}

\newcommand{\m}{\mathfrak m}

\newcommand{\ZZ}{\mathbb{Z}}
\newcommand{\NN}{\mathbb{N}}

\newcommand{\K}{\mathbb{K}}

\newtheorem{lemma}{Lemma}[section]
\newtheorem{corollary}[lemma]{Corollary}
\newtheorem{theorem}[lemma]{Theorem}
\newtheorem{proposition}[lemma]{Proposition}
\newtheorem{definition}[lemma]{Definition}
\newtheorem{remark}[lemma]{Remark}
\newtheorem{example}[lemma]{Example}
\newtheorem{question}[lemma]{Question}

\advance\headheight1.15pt

\begin{document}

\pagenumbering{arabic}
	
	\title[Regularity of symbolic powers of weighted oriented graphs]{Regularity of symbolic and ordinary powers of weighted oriented graphs and their upper bounds} 
	
\author[Manohar Kumar]{Manohar Kumar$^*$}
\address{Department of Mathematics, Indian Institute of Technology
		Kharagpur, West Bengal, INDIA - 721302.}
\email{manhar349@gmail.com}
	
\author[Ramakrishna Nanduri]{Ramakrishna Nanduri$^{\dagger}$}
\address{Department of Mathematics, Indian Institute of Technology
		Kharagpur, West Bengal, INDIA - 721302.}
		
	\thanks{$^*$ Supported by PMRF fellowship, India}
	\thanks{$^\dagger$ Corresponding author and supported by SERB grant No: CRG/2021/000465, India}
	\thanks{AMS Classification 2020: 13D02, 05E40, 05E99, 13D45}
	\email{nanduri@maths.iitkgp.ac.in}
	
	\maketitle

\begin{abstract}
In this paper, we compare the regularities of symbolic and ordinary powers of edge ideals of weighted oriented graphs. For any weighted oriented complete graph $K_n$, we show that $\reg(I(K_n)^{(k)})\leq \reg(I(K_n)^k)$ for all $k\geq 1$. Also, we give explicit formulas for $\reg(I(K_n)^{(k)})$ and $\reg(I(K_n)^{k})$, for any $k\geq 1$. As a consequence, we show that $\reg(I(K_n)^{(k)})$ is eventually a linear function of $k$. For any weighted oriented graph $D$, if $V^+$ are sink vertices, then we show that $\reg(I(D)^{(k)}) \leq \reg(I(D)^k)$ with $k=2,3$ and equality cases studied. Furthermore, we give formula for $\reg(I(D)^2)$ in terms of $\reg(I(D)^{(2)})$ and regularity of certain induced subgraphs of $D$. 

Finally,  we compare the regularity of symbolic powers of weighted oriented graphs $D$ and $D'$, where $D'$ is obtained from $D$ by adding a pendant.  	
\end{abstract}

\section{Introduction}
The study of regularity of powers of homogeneous ideals in a standard graded ring has been a central research focus for decades. A celebrated result, established independently by Cutkosky, Herzog, Trung, and Kodiyalam, showed that for a homogeneous ideal $I$ in a standard graded algebra, regularity function $\reg(I^k)$ is eventually a linear function of $k$, see \cite{CHT99, K00}, whereas regularity of symbolic powers of $I$,  $\reg(I^{(k)})$ need not be a linear function of $k$, for $k \gg 0$ (see \cite[Proposition 7]{ctv93}, \cite[Example 4.4]{CHT99} or \cite[Theorem 5.15]{dhnt21}). However, in \cite{hht07}, Herzog, Hibi, and Trung showed that $\reg(I^{(k)})$ is eventually a quasi-linear function if $I$ is a monomial ideal in a polynomial ring over a field. Whereas, it is not known if the regularity of symbolic
powers of edge ideals of simple graphs is a linear function rather than a quasi-linear
one.   

In this work, we compare regularity of symbolic and ordinary powers of edge ideals of weighted oriented graphs. Weighted oriented graphs play an important role in diverse research areas and various applications in Coding Theory, Commutative Algebra, and Algebraic Geometry etc, see \cite{hlmrv19,prt19}.  

A {\em weighted oriented graph} is a triplet $D=(V(D), E(D), w)$, where $V(D)$ is the vertex set of $D$, $E(D)=\{(x,y) |\mbox{ there is an edge from vertex $x$ to vertex $y$}\}$ is the {\it edge set} of $D$, $w:V(D)\rightarrow \NN$ is a weight function on the vertex set of $D$, and $D$ has no loop edges and no parallel edges. We say that a vertex $x$ of $D$ has nontrivial weight if $w(x)>1$. Let $V(D)=\{x_1,\ldots,x_n\}$ and $R=\mathbb{K}[x_1,\ldots,x_n]$ be the polynomial ring in $n$ variables $x_1,\ldots,x_n$, over a field $\mathbb{K}$. Then the {\it edge ideal} of $D$ is defined as the ideal 
\begin{equation*}
    I(D)=(x_ix_j^{w(x_j)}~|~(x_i,x_j)\in E(D))\subset R.
\end{equation*}
   Many authors are studying regularity of symbolic and ordinary powers of $I(D)$ and giving linear upper bounds for them by fixing certain orientations or weights of $D$. Equality of symbolic and ordinary powers of edge ideals of weighted oriented graphs is studied in \cite{mp21},\cite{gmv21},\cite{mp23},\cite{b23},\cite{k23}. In \cite{mp21}, the authors showed that if $D$ is a weighted oriented odd cycle, then $\reg(I(D)^{(k)}) \leq \reg(I(D)^k)$, for all $k\geq 1$, and equality holds if $D$ has only one vertex with non-trivial weight. In \cite{dhnt21}, Dung, Hien, Nguyen, and Trung gave a bound for the regularity of symbolic powers of monomial ideals in terms of the maximal degree of a vertex of the symbolic polyhedral cone. Later an upper bound for $\reg(I(K_n)^k)$ is given in \cite{czw22} for a fixed orientation on a weighted oriented complete graph $K_n$. In this work, for any weighted oriented graph $D$, we compare the regularities $\reg(I(D)^{(k)}), \reg(I(D)^k)$ for $k\geq 1$.
Furthermore, we have shown that
         $$\reg(I(D)^{(k)}) \leq \reg(I(D')^{(k)}) \mbox{ for all } k \geq 1 , $$
where $D'$ is a weighted oriented graph obtained from any weighted oriented graph $D$ by adding a pendant (Theorem \ref{thmsymboliccomaprision}).         
In this paper, for any weighted oriented graph $D$, under the assumption that $V^+$ are sink vertices, we showed that $\reg(I(D)^{(k)}) \leq \reg(I(D)^k)$, for $k = 2, 3$ (Theorem \ref{thm1}). Furthermore, if $D$ is a weighted oriented graph such that $V^+$ are sink vertices, then we give a formula for the regularity of $I(D)^2$ in terms of regularity of $\reg I(D)^{(2)}$ and regularity of certain induced subgraphs of $D$ (see Theorem \ref{thm3.5} for more detail). This shows that the equality of regularities of symbolic and ordinary powers of $I(D)$ depends on the weights of the vertices. 
For any weighted oriented complete graphs $K_n$, we give explicit formulas for $\reg(I(K_n)^{k})$ and $\reg(I(K_n)^{(k)})$ in terms of the weights of the vertices and also showed that 
$$\reg(I(K_n)^{(k)}) \leq \reg(I(K_n)^k) \mbox{ for all } k \geq 1 ~~(\mbox{Theorem } \ref{thmcompletesymbolic1} \mbox{ and Theorem } \ref{thmcompletesymbolic2}).$$ 
Also, we deduce that $\reg(I(K_n)^{(k)})$ is an eventually linear function of $k$ (Corollary \ref{eventually-linear}).
Our approach uses the regularity formula of monomial ideals in terms of the reduced homology of their degree complexes. 
\vskip 0.2cm 
\noindent
Now, we give section-wise description of the paper. In Section \ref{sec2}, we recall the definitions and essential results required to prove our main results. In Section \ref{sec3}, we compare regularity of symbolic powers of weighted oriented graphs $D$ and $D'$, where $D'$ is obtained from $D$ by adding a pendant. In Section \ref{sec4}, we compare the regularities of small symbolic and ordinary powers of edge ideals of $D$. Finally in Section \ref{sec5}, we give explicit formulas for regularity of symbolic and ordinary powers for any weighted oriented complete graph $K_n$.
	
\section{Preliminaries} \label{sec2}

  In this section, we recall some definitions and results which will be used throughout the paper.  \\ 
 Let $D=(V(D), E(D), w)$ be a weighted oriented graph and $V(D)=\{x_1,\ldots,x_n\}$. Then the  underlying simple graph $G$ of $D$ is a simple graph such that $V(G)=V(D)$ and $E(G)=\{\{x,y\}| (x,y) \mbox { or } (y,x)\in E(D)\}$. That is, $G$ is the simple graph without orientation and weights in $D$. Let $R=\mathbb{K}[x_1,\ldots,x_n]$, and $\mathfrak{m}=(x_1,\ldots,x_n)$, the homogeneous maximal ideal in $R$, where $\mathbb{K}$ is a field. For a vertex $x\in V(D)$, its {\it outer neighbourhood} is defined as 
  $N_D(x)^+ := \{y \in V(D) | (x,y)\in E(D)\}$, its {\it inner neighbourhood} is defined as $N_D(x)^- := \{z \in V(D) | (z,x)\in E(D)\}$, and $N_D(x):=N_D(x)^+ \cup N_D(x)^-$. Also, $N_D[x]:=N_D(x)^+ \cup N_D(x)^- \cup \{x\}$, $N_D[x]^+:=N_D(x)^+\cup \{x\}$, $N_D[x]^-:=N_D(x)^-\cup \{x\}$. A vertex $x\in V(D)$ is called a {\it source} if $N_D(x)^-=\emptyset$ and $x$ is called a {\it sink} if $ N_D(x)^+ = \emptyset $. If $x\in V(D)$ is a source, then we set $w(x)=1$. Let $V^+(D):= \{x\in V(D): w(x)\geq 2\}$, simply denoted by $V^+$. 
  The degree of a vertex $x\in V(D)$ is defined as $\text{d}_D(x) :=|N_D(x)|$.
  A vertex $x$ is called a {\it leaf } if its degree is one. A set of vertices $C$ of $G$ is called a {\it vertex cover} of $G$ if $C\cap e\neq \emptyset$ for all $e\in E(G)$. For a vertex cover $C$ of a weighted oriented graph $D$, define 
$L_1^D(C):= \{x \in C: N_D^+(x)\cap C^c \neq \emptyset\}$,
$L_2^D(C):= \{x \in C: x \notin L_1^D(C), N_D^-(x)\cap C^c \neq \emptyset\}$,
$L_3^D(C):= C\setminus (L_1^D(C)\cup L_2^D(C))$, where $C^c=V(D)\setminus C$.
A vertex cover $C$ of a weighted oriented graph $D$ is called {\it strong}, if for each $x\in L_3^D(C)$, there exists $y \in (L_2^D(C)\cup L_3^D(C))\cap V^+$ such that $(y,x)\in E(D)$. The {\it irreducible ideal associated to} a vertex cover $C$ of $D$ is the ideal
\begin{equation*}
    I_C^D=(L_1^D(C)\cup \{x_j^{w(x_j)} \mid x_j \in L_2^D(C) \cup L_3^D(C)\}) \subseteq R.
\end{equation*} 
 Let $\Delta$ be a simplicial complex on $[n]=\{1,\ldots,n\}$ that is a collection of subsets of $[n]$ closed under taking subsets. The dimension of $\Delta$ is defined as $\dim \Delta := \max\{\dim F \mid F \in \Delta\}$, where $\dim F = \vert F \rvert -1 $ and $\vert F \rvert$ is the cardinality of $F$. The set of its maximal elements under inclusion, called by facets, is denoted by $\mathfrak{F}(\Delta)$. For any $S\subseteq V(D)$, we denote $D\setminus S$, the induced subgraph of $D$ on the vertex set $ V(D)\setminus S $. For a face $F \in \Delta$, the link of $F$ in $\Delta$ is the subsimplicial complex of $\Delta$ defined by 
  \begin{equation*}
      \text{lk}_{\Delta}F := \{G \in \Delta \mid F\cup G \in \Delta, F \cap G=\emptyset\}.
  \end{equation*} 
\noindent
Now, we recall degree complex of a monomial ideal with respect to $\mathbf{a} \in \ZZ^n$ which was used by Takayama in \cite{t05} to give a combinatorial formula for $\dim
_{K}(H_{\mathfrak{m}}^i(R/I)_{\mathbf{a}})$, for a monomial ideal $I$ in $R$, where $H_{\mathfrak{m}}^i(R/I)$ is the $i^{th}$ local cohomology. 
For $\mathbf{a}=(a_1,\ldots,a_n) \in \ZZ^n$, set ${\bf x}^{\bf a}=x_1^{a_1}\cdots x_n^{a_n}$. We denote by $G_{\mathbf{a}}$ the negative support of $\mathbf{a}$, that is  $G_{\mathbf{a}} = \{i \in [n] \mid a_i < 0\}$. For every subset $F \subseteq [n]$, let $R_F=R[x_i^{-1} \mid i \in F]$. The {\it degree complex} of $I$ with respect to $\mathbf{a}$ is defined as
\begin{equation*}
    \Delta_{\mathbf{a}}(I) := \{F \setminus G_{\mathbf{a}} \mid G_{\mathbf{a}} \subseteq F, {\bf x}^{\bf a} \not \in IR_F \}.
\end{equation*}
 The $\Delta_{\mathbf{a}}(I)$ is a simplicial complex on the vertex set $[n]$.
For a monomial $f \in R=\mathbb{K}[x_1,\ldots,x_n]$,  define its {\it support} as $\supp(f) := \{x_i : x_i \mid f\}$. We denote $\mathcal{G}(I)=$ minimal set of monomial generators of $I$. For a monomial ideal $I$, its support is defined as $\supp(I):= \displaystyle \cup_{f\in \mathcal{G}(I)} \supp(f)$. For an exponent $\mathbf{a} \in \ZZ^n$, we denote $\supp(\mathbf{a})=\{i \in [n] \mid a_i \neq 0\}$, the support of $\mathbf{a}$. For any subset $U \subseteq [n]$, we denote 
\begin{equation*}
    I_U :=(f \mid f \text{ is a monomial which belongs to } I \text{ and } \supp(f) \subseteq U)
\end{equation*}
an ideal in $R$ known as the restriction of $I$ to $U$. 

 \begin{definition} For any homogeneous ideal $I$ in $R$, define the Castelnuovo-Mumford regularity (or merely regularity) as 
 \begin{eqnarray*}
 \reg(I) &=&  \max \{j - i \mid \beta_{i,j}(I) \neq 0\} \\
         &=& \max\{j+i \mid H_{\m}^i(I)_j \neq 0\},   
 \end{eqnarray*}
  where $\beta_{i,j}(I)$ is the $(i,j)^{th}$ graded Betti number of $I$ and $H_{\m}^i(I)_j$ denotes the $j^{th}$ graded component of the $i^{t h}$ local cohomology module $H_{\m}^i(I)$.  
\end{definition}

 \begin{definition}[Symbolic power] Let $R$ be a Noetherian ring and $I\subset R$ be an ideal. Then the  $k$-th symbolic power of $I
 $ is defined as $$ I^{(k)}=\bigcap_{P\in Ass(R/I)}(I^k R_P \cap R).$$
 \end{definition}   

\begin{definition}[Extremal exponent] Let $I$ be a monomial ideal in $R$. A pair $({\bf a},i)\in \NN^{n}\times \NN$ is called an extremal exponent of $I$, if there exists a face $F \in \Delta_{{\bf a}}(I)$ with $F \cap \supp {\bf a} = \emptyset $ such that $\Tilde{H}_{i-1}(\mathrm{lk}_{\Delta_{\mathbf{a}}(I)}F; \K)\neq 0 $ and $\reg(R/I)=\vert {\bf a} \vert +i $, where $\Tilde{H}_{i}(-;\K)$ denotes the $i^{th}$ reduced homology over $\K$. 
\end{definition}
\noindent 
Below, we recall some lemmas which we use frequently in many proofs. 

\begin{lemma}\cite[Lemma 2.12]{m22}\label{hom1}
 Let $I$ be a monomial ideal in $R$. Then,
 \begin{align*}
   \reg(R/I)= & \max\{|\mathbf{a}|+i \mid \mathbf{a} \in \mathbb{N}^n, i\geq 0, \Tilde{H}_{i-1}(\mathrm{lk}_{\Delta_{\mathbf{a}}(I)}F; \K)\neq 0 
    \mbox{ for some } \\ 
     & F \in \Delta_{\mathbf{a}}(I) \mbox{ with } F \cap \supp \mathbf{a} = \emptyset \}. 
 \end{align*}
 In particular, if $I=I_{\Delta}$ is the Stanley-Reisner ideal of a simplicial complex $\Delta$, then
 \begin{equation*}
  \reg(\K[\Delta])=\max\{i \mid i\geq 0, \Tilde{H}_{i-1}(\Delta_{\mathbf{a}}(I); \K)\neq 0   \mbox{ for some } F \in \Delta\}.   
 \end{equation*} 
\end{lemma}

\begin{lemma}\cite[Lemma 2.19]{m22}\label{lm2.1}
 Let $I, J$ be two proper monomial ideals of $R$. Let $(\mathbf{a},i)$ be an extremal exponent of $I$. If $\Delta_{\mathbf{a}}(I)=\Delta_{\mathbf{a}}(J)$, then $\reg(I) \leq \reg(J)$. In particular, if $J \subseteq I$ and $\Delta_{\mathbf{a}}(I)=\Delta_{\mathbf{a}}(J)$ for all exponent $\mathbf{a} \in \mathbb{N}^n$ such that ${\bf x}^{\bf a} \not \in I$, then $\reg(I) \leq \reg(J)$. 
\end{lemma}
\begin{lemma}\cite[Lemma 2.18]{m22}\label{lm2.2}
 Let $I$ be a monomial ideal in $R=\mathbb{K}[x_1,\ldots,x_n]$ and $\mathbf{a} \in \mathbb{N}^n$. Then,
 \begin{equation*}
     I_{\Delta_\mathbf{a}(I)}=\sqrt{I:{\bf x}^{\bf a}}.
 \end{equation*}
 In particular, ${\bf x}^{\bf a} \in I$ if and only if $\Delta_{\mathbf{a}}(I)$ is the void complex.
\end{lemma}
\begin{lemma}\cite[Lemma 2.23]{m22}\label{lm2.3}
Let $I$ be a monomial ideal and $(\mathbf{a},i) \in \mathbb{N}^n \times \mathbb{N}$ be an extremal exponent of $I$. Assume that $x_j$ is a variable that appears in $\sqrt{I:{\bf x}^{\bf a}}$ and $j \not \in \supp(\mathbf{a})$. Then,
\begin{equation*}
    \reg(I)=\reg(I,x_j)=\reg(I_U),
\end{equation*}
where $U=[n]\setminus \{j\}$.
\end{lemma}

\begin{lemma}\cite[Corollary 2.22]{m22}\label{lm2.4}
  Let $J$ be a monomial ideal in $R$. Let $U \subseteq [n]$. Then,
  \begin{equation*}
      \reg(J_U) \leq \reg(J). 
  \end{equation*}
\end{lemma}

\begin{lemma}\cite[Lemma 2.24]{m22}\label{lm2.6}
 Let $I$ be a monomial ideal in $R=\mathbb{K}[x_1,\ldots,x_n]$ generated by the monomial $f_1,\ldots,f_r$ and $\mathbf{a} \in \mathbb{N}^n$. Then,  
$\sqrt{I : {\bf x}^{\bf a}} $ is generated by $\sqrt{f_1/\gcd(f_1,{\bf x}^{\bf a})}, \ldots, \\
\sqrt{f_r/\gcd(f_r,{\bf x}^{\bf a})}$.
\end{lemma} 

\begin{lemma}\cite[Lemma 2.26]{m22}\label{lm2.7}
 Let $J \subseteq L$ be two monomial ideals in $R$ and $\mathbf{a}$ be an exponent such that ${\bf x}^{\bf a} \not \in L$. Let $g$ be a minimal generator of $L$ and $f=\sqrt{g/\gcd(g,{\bf x}^{\bf a})}$. Assume that $f$ belongs to $\sqrt{J :g}$. Then, $f$ belongs to $\sqrt{J:{\bf x}^{\bf a}}$.
\end{lemma}

\begin{lemma}\cite[Corollary 3.12]{s08}\label{lm2.9}
For any simple graph $G$, 
\begin{equation*}
    I(G)^{(2)}=I(G)^{\{2\}}+I(G)^2, 
\end{equation*}
where $I(G)^{\{2\}}$ is known as $2^{nd}$ secant ideal.  
In particular, $I(G)^{(2)}$ is generated by cubics of the form $x_ix_jx_k$ such that $\{x_i,x_j,x_k\}$ is a triangle in $G$ and quadrics of the form $x_ix_jx_kx_l$ such that $\{x_i,x_j\}$ and $\{x_k,x_l\}$ are edges of $G$.
\end{lemma}

Recall the following lemma, known as regularity lemma, which we frequently use in many proofs. 
\begin{lemma}\cite[Lemma 1.2]{htt16} [Regularity lemma] \label{lm2.12}
		Let $0\rightarrow A \rightarrow B \rightarrow C \rightarrow 0$ be a short exact sequence of finitely generated graded $R$-modules. Then
		\begin{enumerate}
			\item $\reg(B) \leq \max\{\reg(A), \reg(C)\}$, 
			\item $\reg(A) \leq \max\{\reg(B), \reg(C)+1\}$,
			\item $\reg(C) \leq \max\{\reg(A)-1, \reg(B)\}$. 

		\end{enumerate}
\end{lemma}
\noindent
  Let $ D $ be a  weighted oriented graph such that the vertices of $V^{+}$ are sinks. Recall from \cite[Notation 3.6]{mp21} that the map 
	\begin{equation*}
		\Phi : R \longrightarrow R, ~ \mbox{defined as }  ~ x_j \mapsto x_j ~ \mbox{if} ~ x_j \notin V^{+},  \mbox{and} ~ x_j \mapsto x_j^{w(x_j)}~  \mbox{if}  ~ x_j \in V^{+}. 
	\end{equation*} 
 
\begin{lemma}\cite[Theorem 3.7]{mp21}\label{lm2.8}
 Let $D$ be a weighted oriented graph where the elements of  $V^+$ are sinks and $G$ be the underlying simple graph of $D$. Then, $\Phi(I(G)^k)= I(D)^k$ and $\Phi(I(G)^{(k)})=I(D)^{(k)}$ for all $k\geq 1$. 
\end{lemma}

\section{Regularity comparison of symbolic powers of weighted oriented graphs}\label{sec3}

In this section, we compare the regularity of symbolic powers of weighted oriented graphs $D$ and $D'$, where $D'$ is obtained from $D$ by adding a pendant. First, we give a remark for a lower bound of $\reg(I(D)^{(k)})$ if $V^+$ are sinks.
\begin{remark}\label{remark-lower}
Let $D$ be a weighted oriented graph such that $V^+$ are sinks and $I=I(D)$ be its edge ideal. Suppose $H$ is an induced subgraph of $G$ and $J=\Phi(I(H))R$.
Note that for edge ideals of a simple graph $G$ and its induced subgraph $H$, we have $I(H)^{(k)}=I(G)^{(k)}\vert_U$, where $U=V(H)$. This implies that $\Phi(I(H)^{(k)})=\Phi(I(G)^{(k)}\vert_U)=\Phi(I(G)^{(k)})\vert_U$ because $\Phi$ preserves the support. Therefore, by Lemma \ref{lm2.8}, we have $J^{(k)}=I^{(k)}\vert_U$. Thus, $J^{(k)}$ is a retraction of $I^{(k)}$. Therefore, by \cite[Corollary 2.5]{ohh00} we have that 
$\beta_{i,j}(J^{(k)}) \leq \beta_{i,j}(I^{(k)})$ for all  $i,j \geq 0$ and $k\geq 1$. In particular, let $D_1$ be an induced matching in $D$, say $E(D_1) = \{(x_1, y_1), \ldots, (x_r, y_r)\}$. Then, we have
 \begin{equation*}
(k-1)(w+1)+\sum_{i=1}^rw(y_i)+1=\reg(I(D_1)^{k})=\reg(I(D_1)^{(k)}) \leq \reg(I(D)^{(k)}), \mbox{ for all } k\geq 1. 
 \end{equation*}
\end{remark}
\noindent 
The Example \ref{example3.6} at the end of this section shows that if we drop the assumption that $V^+$ are sinks, then the inequality in Remark \ref{remark-lower} may not be true.
Now, we prove that $\reg(I(D)^{(k)}) \leq \reg(I(D')^{(k)})$ for all $k\geq 1$, where $D'$ is obtained from $D$ by adding a pendant. To prove that first we prove two lemmas and a proposition.
 
\begin{lemma}\label{lemstrongcover2}
Let $D$ be a weighted oriented graph having at least one sink vertex $x$ with $w(x) \geq 2$. Suppose $D'$ is a weighted oriented graph obtained by adding an edge $(x,y)$ to $D$, where $y$ is a new vertex which is a leaf. If $C$ is a strong vertex cover of $D$, then $C\cup \{y\}$ is a strong vertex cover of $D'$. Moreover, a strong vertex cover of $D'$ containing the vertex $y$ is of the form $C \cup \{y\}$ for some strong vertex cover $C$ of $D$. 
\end{lemma}
\begin{proof} Let $C$ be a strong vertex cover of $D$.  If $x \not \in C$, then $y \in L_2^{D'}(C\cup \{y\})$ because $N^-_{D'}(y)=\{x\} \not \subseteq C\cup \{y\}$. Thus, $C \cup \{y\}$ is a strong vertex cover of $D'$. Assume $x \in C$. Since $x$ is a sink vertex of $D$, we have $x \in L_2^{D}(C) \cup L_3^{D}(C)$ and $y \in L_3^{D'}(C\cup \{y\})$. Thus, $C \cup \{y\}$ is a strong vertex cover of $D'$ because  $(x,y) \in E(D')$ with $x \in (L_2^{D'}(C\cup \{y\}) \cup L_3^{D'}(C \cup \{y\})) \cap V^+$.  \\
Now, suppose $C'$ is a strong vertex cover of $D'$ containing $y$. Then we have to show that $C' \setminus y$ is a strong vertex cover of $D$. If $x \not\in C'$, then $C' \setminus \{y\}$ is a vertex cover of $D$ which is strong, as required. Assume $x \in C'$. Then, $x \in L_2^{D}(C'\setminus y)\cup L_3^{D}(C'\setminus y)$ because $x$ is a sink vertex of $D$. If $ x \in L_2^{D}(C' \setminus \{y\})$, then $C'\setminus y$ is a strong vertex cover of $D$. Now, suppose  $ x \in L_3^{D}(C' \setminus \{y\})$. Then, by \cite[Proposition 5]{prt19}, $N_D(x) \subseteq C'\setminus y$. This implies that $x \in L_3^{D'}(C')$ because $N_{D'}(x) \subseteq C'$. Since $C'$ is a strong vertex cover of $D'$, this implies that there exists an edge $(z,x) \in E(D')$ with $z\neq y$ such that $z \in L_2^{D'}(C')\cup L_3^{D'}(C') \cap V^+$. This gives that $(z,x) \in E(D)$ and $z \in L_2^{D}(C'\setminus \{y\})\cup L_3^{D}(C'\setminus \{y\}) \cap V^+$. Thus, $C'\setminus y$ is a strong vertex cover of $D$.
  \end{proof}
  \begin{lemma}\label{lemintersectionsum}
Let $J$ and $I_i$  for $i \in S=\{1 \ldots u \} $ be any monomial ideal such that $\supp{(J)} \cap \supp{(I_i)} = \emptyset$ for all $i=1\ldots u$. Then,  
\begin{align*}
    \bigcap_{i \in S}(J+I_i)^k =\sum_{a=0}^kJ^{a}\bigcap_{i \in S}I_i^{k-a}, \text{ for all $k\geq 1$. }
\end{align*} 
\end{lemma}
\begin{proof}
Using binomial expansion, we have
\begin{align*}
 \bigcap_{i \in S}(J+I_i)^k &=\bigcap_{i \in S} \sum_{k_1+k_2=k}J^{k_1}I_i^{k_2} \\
   &=\sum_{a_1+b_1=k}J^{a_1}I_1^{b_1} \cap \ldots \cap \sum_{a_u+b_u=k}J^{a_u}I_u^{b_u}\\
 &=\sum_{a_1+b_1=k}(J^{a_1} \cap I_1^{b_1}) \cap \ldots \cap \sum_{a_u+b_u=k}(J^{a_u} \cap I_u^{b_u})  \text{ (because  
                      }
\supp{(J)} \cap \supp{(I_i)} = \emptyset)\\
 &=\sum_{a_1+b_1=k} \ldots \sum_{a_u+b_u=k} (J^{a_1} \cap I_1^{b_1}) \cap \ldots \cap (J^{a_u} \cap I_u^{b_u}) \\
\nonumber        &=\sum_{a_1+b_1=k} \ldots \sum_{a_u+b_u=k} (J^{a_1} \cap \ldots \cap J^{a_u}) \cap (I_1^{b_1} \cap \ldots \cap I_u^{b_u}) \\
\nonumber  &=\sum_{a_1+b_1=k} \ldots \sum_{a_u+b_u=k} (J^{a_1} \cap \ldots \cap J^{a_u})(I_1^{b_1} \cap \ldots \cap I_u^{b_u})\\
& \text{ (because  
                      }
\supp{(J)} \cap \supp{(I_i)} = \emptyset) \\
&=\sum_{a_1+b_1=k} \ldots \sum_{a_u+b_u=k} J^{\max\{a_1,\ldots,a_u\}}(I_1^{b_1} \cap \ldots \cap I_u^{b_u}) \\
& \subseteq \sum_{a_1+b_1=k} \ldots \sum_{a_u+b_u=k} J^a(I_1^{b_1-a} \cap \ldots \cap I_u^{b_u-a}) \text{ (where
                      } a=\max\{a_1,\ldots,a_u\}) \\
 &=\sum_{a=0}^kJ^{a}(I_1^{k-a} \cap \ldots \cap I_u^{k-a})
\end{align*} 
Note that the reverse containment always holds (by taking $a_1=a_2=\cdots=a_u$). 

Therefore, we have
\begin{align*}
    \bigcap_{i \in S}(J+I_i)^k =\sum_{a=0}^kJ^{a}\bigcap_{i \in S}I_i^{k-a}, \text{ for all $k\geq 1$. }
\end{align*} 
\end{proof}
 
 \begin{proposition}\label{proprestriction}
Let $D$ be a weighted oriented graph having at least one sink vertex $x$ with $w(x) \geq 2$. Suppose $D'$ is a weighted oriented graph obtained by adding an edge $(x,y)$ to $D$, where $y$ is a new vertex which is a leaf. Then, $I(D')^{(k)}=I(D)^{(k)}+y^{w(y)}J$, for some monomial ideal $J\subseteq R $. 
 \end{proposition} 
\begin{proof}
Note that by \cite[Theorem 3.7]{prt19}, $I(D)^{(k)}=\displaystyle\cap_{C\in \mathcal{C}_s} I_C^k$, where $\mathcal{C}_s$ is a set of all strong vertex cover $C$ of $D$. Now, suppose  $\mathcal{C'}_s$ is a set of all strong vertex cover $C'$ of $D'$ and $\mathcal{C'}_{s_y}$ is a set of all strong vertex cover $C'$ not containing $y$ of $D'$.  We can write,
 \begin{align*}
   I(D')^{(k)}&=\displaystyle \bigcap_{C'\in \mathcal{C'}_s} (I_{C'}^{D'})^k  \\
              &=\displaystyle \bigcap_{C\cup \{y\} \in \mathcal{C'}_s} (I_{C'}^{D'})^k \cap \left(\bigcap_{C' \in \mathcal{C'}_{s_y} } (I_{C'}^{D'})^k \right) \mbox{ (by Lemma \ref{lemstrongcover2})} \\
              &=\displaystyle \bigcap_{C\in \mathcal{C}_s} (y^{w(y)}+I_C^D)^k \cap \left(\bigcap_{C' \in \mathcal{C'}_{s_y }} (I_{C'}^{D'})^k \right)  \\
              &= \displaystyle  \sum_{a=0}^k\left(y^{w(y)a} \bigcap_{C\in \mathcal{C}_s}(I_C^D)^{k-a}\right)\cap \left(\bigcap_{C' \in \mathcal{C'}_{s_y} } (I_{C'}^{D'})^k \right)  \mbox{ (by Lemma \ref{lemintersectionsum})}  \\
               &= \displaystyle  \sum_{a=0}^k\left(y^{w(y)a} I(D)^{(k-a)}\right)\cap \left(\bigcap_{C' \in \mathcal{C'}_{s_y} } (I_{C'}^{D'})^k \right) \\
           &= \displaystyle       I(D)^{(k)} \cap \left(\bigcap_{C' \in \mathcal{C'}_{s_y} } (I_{C'}^{D'})^k \right)+y^{w(y)} \sum_{a=0}^{k-1}(y^{w(y)a} I(D)^{(k-1-a)})\cap \left(\bigcap_{C' \in \mathcal{C'}_{s_y} } (I_{C'}^{D'})^k \right).
 \end{align*}
Note that for $C' \in \mathcal{C'}_{s_y}$, we have $x \in C'$ because $y \not \in C'$. This implies that $C'$ is a strong vertex cover of $D$ by definition. Thus, $\mathcal{C'}_{s_y} \subseteq \mathcal{C}_s$. Thus, for any $C' \in \mathcal{C'}_{s_y}$,  $x \in L_1^{D'}(C')$ whereas $x \in L_2^{D}(C') \cup L_3^{D}(C')$ because $x$ is sink in $D$ but not in $D'$. Then for any $C' \in \mathcal{C'}_{s_y}$, $I_{C'}^D \subseteq I_{C'}^{D'}$, where $I_{C'}^D$ and $I_{C'}^{D'}$ are irreducible ideals associated to $C'$ in $D$ and $D'$ respectively. Therefore, we can write,
\begin{align*}
I(D)^{(k)}=\cap_{C\in \mathcal{C}_s} (I_C^D)^k 
          = \bigcap_{C' \in \mathcal{C'}_{s_y} } (I_{C'}^D)^k \cap \left(\bigcap_{C \in \mathcal{C}_s\setminus \mathcal{C'}_{s_y} } (I_{C}^D)^k \right).  
\end{align*}
This implies that
$I(D)^{(k)} \subseteq \bigcap_{C' \in \mathcal{C'}_{s_y} } (I_{C'}^D)^k  \subseteq \bigcap_{C' \in \mathcal{C'}_{s_y} } (I_{C'}^{D'})^k.$ \\
Thus, $I(D')^{(k)}=I(D)^{(k)}+y^{w(y)}J$, for some monomial ideal $J\subseteq R$, as required. 
\end{proof}
\noindent 
Now, we prove the main result of this section.
\begin{theorem}\label{thmsymboliccomaprision}
 Let $D$ and $D'$ be any weighted oriented graph as in Proposition \ref{proprestriction}. Then,
 \begin{equation*}
      \reg(I(D)^{(k)}) \leq \reg(I(D')^{(k)}), \mbox{  for all 
  }k\geq 1.
 \end{equation*}
\end{theorem}
\begin{proof}
  By Proposition \ref{proprestriction}, we have $(I(D)^{(k)},y)=(I(D')^{(k)},y)$. Thus, by \cite[Lemma 2.7(1)]{x21} and \cite[Corollary 4.8]{chhktt19}, we have 
  \begin{equation*}
      \reg(I(D)^{(k)})=  \reg((I(D)^{(k)},y)) =\reg((I(D')^{(k)},y)) \leq \reg(I(D')^{(k)})  \mbox{ for all } k \geq 1.
 \end{equation*}
\end{proof} 

The following example illustrates that if we add a new vertex $y$ to $D$ as a leaf, but $y$ is not a sink, then Theorem \ref{thmsymboliccomaprision} may not hold.
\begin{example}\label{example3.6}
	Let $D$ be a Cohen-Macaulay weighted oriented forest as in Figure \ref{fig1}.
\begin{figure}[!h]
    \centering
     \includegraphics[width=0.2 \textwidth]{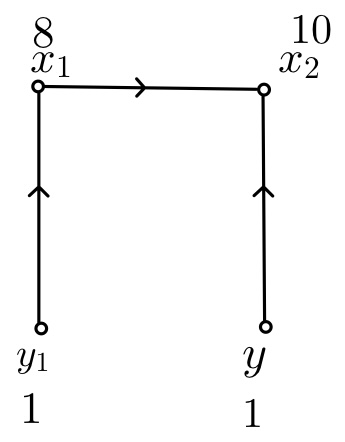} 
    \caption{Cohen-Macaulay weighted oriented forest}
    \label{fig1}
\end{figure}
\\
Then the edge ideal $I(D)=(x_1x_2^{10},y_1x_1^8,yx_2^{10}).$ Also, $D\setminus y$ 
 is an induced subgraph with edge ideal $I(D\setminus y)=(x_1x_2^{10},y_1x_1^8)$. 
Using Macaulay2 \cite{gs}, we have 
\begin{align*}
29=\reg(I(D\setminus y)^{(2)}) &> \reg(I(D)^{(2)})=28 \\
40= \reg(I(D\setminus y)^{(3)})  &> \reg(I(D)^{(3)})=38.   
\end{align*}
\end{example}

\section{Comparing regularity of small Symbolic and ordinary powers}\label{sec4}
In this section, for any weighted oriented graph $D$, we show that $\reg(I(D)^{(k)}) \leq \reg(I(D)^k)$ for $k=2,3$ by assuming that $V^+$ are sink vertices.

\begin{lemma}\label{lm4.1}
Let $I=I(D)$ and $\mathbf{a} \in \mathbb{N}^n$ such that ${\bf x}^{\bf a} \not \in I^{(2)}$. Assume that $V^+$ are sink vertices. Then, $\sqrt{I^{(2)}:{\bf x}^{\bf a}} = \sqrt{I^{2}:{\bf x}^{\bf a}}$.
\end{lemma}
\begin{proof} Let $G$ be the underlying simple graph of $D$. 
Note that $I^2 \subseteq I^{(2)}$. This implies that $\sqrt{I^{2}:{\bf x}^{\bf a}} \subseteq \sqrt{I^{(2)}:{\bf x}^{\bf a}}$. Now, we show the reverse inclusion. Let $f \in \mathcal{G}(\sqrt{I^{(2)} : {\bf x}^{\bf a}} )$. By Lemma \ref{lm2.9}, we have that  $I(G)^{(2)}=I(G)^2+(x_ix_jx_r~|~ \{x_i,x_j,x_r\} \mbox{ is a triangle in  } G)$. Thus, by Lemma \ref{lm2.8}, we have
\begin{eqnarray*}
  I^{(2)} &=& \Phi(I(G)^{(2)}) \\
          &=& \Phi(I(G)^2)+\Phi(x_ix_jx_r~|~ \{x_i,x_j,x_r\} \mbox{ is a triangle in } G) \\  
          &=& I^2 +(x_i^{w(x_i)}x_j^{w(x_j)}x_r^{w(x_r)}~|~ \{x_i,x_j,x_r\} \mbox{ is a triangle in } G).              
\end{eqnarray*}
Now, by Lemma \ref{lm2.6} we assume that $f=\sqrt{x_i^{w(x_i)}x_j^{w(x_j)}x_r^{w(x_r)}/\gcd(x_i^{w(x_i)}x_j^{w(x_j)}x_r^{w(x_r)},{\bf x}^{\bf a})}$, where $\{x_i,x_j,x_r\}$ is a triangle in $G$. Note that $x_l^{w(x_l)}x_i^{w(x_i)}x_j^{w(x_j)}x_r^{w(x_r)} \in I^2$ for all $l \in \supp(x_i^{w(x_i)}x_j^{w(x_j)}x_r^{w(x_r)})$. Thus, by Lemma \ref{lm2.7}, we have $f \in \mathcal{G}(\sqrt{I^{2} : {\bf x}^{\bf a}}) $, as required. 
\end{proof}

\begin{lemma}\label{lm4.2}
Let $I=I(D)$ and $\mathbf{a} \in \mathbb{N}^n$ such that ${\bf x}^{\bf a} \not \in I^{(3)}$. Assume that $V^+$ are sink vertices and $\sqrt{I^{(3)}:{\bf x}^{\bf a}} \neq \sqrt{I^{3}:{\bf x}^{\bf a}}$. Let $f\in \mathcal{G}(\sqrt{I^{(3)}:{\bf x}^{\bf a}})$ and $f\notin \sqrt{I^{3}:{\bf x}^{\bf a}}$. Then,
\begin{enumerate}
    \item there exists a triangle $x_ix_jx_r$ of underlying graph $G$ and $v \not \in \{x_i,x_j,x_r\}$ such that $x_i^{w(x_i)}x_j^{w(x_j)}x_r^{w(x_r)}x_l^{w(x_l)} \mid {\bf x}^{\bf a} $,
 \item $\deg f = 1 $ and $\supp f \not \subseteq \supp \mathbf{a}$.
\end{enumerate}
\end{lemma}
\begin{proof}
 By \cite[Theorem 2.10]{m22} we have that $I(G)^{(3)}=I(G)^3 + IJ_1(G)+J_2(G)$, where $J_1(G)$ is the ideal generated by all squarefree monomials $x_ix_jx_r$ such that  $\{x_i,x_j,x_r\}$ forms a triangle, and $J_2(G)$ is the ideal generated by all squarefree monomials $x_sx_tx_ux_l$ such that $\{x_s,x_t,x_u,x_l\}$ forms a clique of size four and all squarefree monomials $x_C$ such that $C$ is a five cycle of $G$. By Lemma \ref{lm2.8}, this implies that 
\begin{equation*}
  I^{(3)}=\Phi(I(G)^{(3)})=\Phi(I(G)^3)+ \Phi(I(G)J_1(G))+\Phi(J_2(G))= I^3+I \widetilde{J_1(G)}+\widetilde{J_2(G)},
\end{equation*}
where $\widetilde{J_1(G)}$ is the ideal generated by monomials $x_i^{w(x_i)}x_j^{w(x_j)}x_r^{w(x_r)}$ such that  $\{x_i,x_j,x_r\}$ forms a triangle, and $\widetilde{J_2(G)}$ is the ideal generated by monomials $x_s^{w(x_s)}x_t^{w(x_t)}x_u^{w(x_u)}x_l^{w(x_l)}$ such that $\{x_s,x_t,x_u,x_l\}$ forms a clique of size four and monomials $x_C^{w(x_C)}$ such that $C$ denotes five cycle of $G$ and $w(x_C)$ denotes its respective weights. 
By Lemma \ref{lm2.6}, there exists a minimal generator $g$ of $I^{(3)}$ such that $f=\sqrt{g/\gcd(g,{\bf x}^{\bf a})}$. Since ${\bf x}^{\bf a} \not \in I^{(3)}$, we have 
\begin{equation}\label{eq3.1}
    \emptyset \neq \supp(f) \subseteq \supp(g) \subseteq \supp(f) \cup \supp \mathbf{a}.
\end{equation}
Now, we consider three possible cases.
\vskip 0.2cm 
\noindent 
{\bf Case 1.} Suppose $g=x_C^{w(x_C)}$, where $C$ denotes a five cycle of $G$ and $w(x_C)$ denotes its respective weights. By equation \eqref{eq3.1} and the fact that $x_j^{w(x_j)}x_C^{w(x_C)} \in I^3$ for any $j \in \supp C$, we have $f \in \sqrt{I^3:x_C^{w(x_C)}}$. By Lemma \ref{lm2.7}, we have $f \in \sqrt{I^3:{\bf x}^{\bf a}}$, a contradiction. Therefore, this case never arises. 
\vskip 0.2cm 
\noindent
{\bf Case 2.} Suppose $g=x_s^{w(x_s)}x_t^{w(x_t)}x_u^{w(x_u)}x_l^{w(x_l)}$, where $\{x_s,x_t,x_u,x_l\}$ forms a clique of size four of $G$. By equation \eqref{eq3.1} and the fact that $x_j^{2w(x_j)}x_s^{w(x_s)}x_t^{w(x_t)}x_u^{w(x_u)}x_l^{w(x_l)}$ for any $j \in \supp(x_s^{w(x_s)}x_t^{w(x_t)}x_u^{w(x_u)}x_l^{w(x_l)}) $, we have $f \in \sqrt{I^3:x_s^{w(x_s)}x_t^{w(x_t)}x_u^{w(x_u)}x_l^{w(x_l)}}$. By Lemma \ref{lm2.7}, $f \in \sqrt{I^3: {\bf x}^{\bf a}}$, a contradiction. Therefore, this case also never arises. 
\vskip 0.2cm 
\noindent
{\bf Case 3.} Suppose $g \in I\widetilde{J_1(G)}$. Then we assume that $g=x_i^{w(x_i)}x_j^{w(x_j)}x_r^{w(x_r)}x_u^{w(x_u)}x_l^{w(x_l)}$, where $x_ix_jx_r \in J_1(G)$  and $x_ux_l \in I(G)$. By Lemma \ref{lm2.7} and the fact that $x_l^{w(x_l)}g \in I^3 $ for any $x_l \in N[{x_i,x_j,x_r}]$, we have $x_l \in \sqrt{I^3:{\bf x}^{\bf a}}$. Since $f \not \in \sqrt{I^3:{\bf x}^{\bf a}} $, we have $x_l \notin \supp(f)$. Therefore, we must have
\begin{equation}\label{eq3.2}
    \supp(f)\cap N[\{x_i,x_j,x_r\}]=\emptyset.
\end{equation}
By equations \eqref{eq3.1}, \eqref{eq3.2}, we have $\supp f \subseteq \{x_u,x_l\} $. If possible, let $f=x_ux_l \in I(G)$. Then, $(x_ux_l)^{\max\{w(x_u),w(x_l)\}} \in I$. This implies that $(x_ux_l)^{\max\{w(x_u),w(x_l)\}}g \in I^3$. Thus, by Lemma \ref{lm2.7},  $f=x_ux_l  \in \sqrt{I^3:{\bf x}^{\bf a}}$ which contradicts our assumption that $f \not \in \sqrt{I^3:{\bf x}^{\bf a}}  $. Therefore, $f \not \in I(G)$. Thus, $\deg f =1$. Without loss of generality, we take $f=x_u$. By equations \eqref{eq3.1}, \eqref{eq3.2}, and the assumption that $f=\sqrt{g/\gcd(g,{\bf x}^{\bf a})}$, we get that 
\begin{align*}
    x_l \not \in \{x_i,x_j,x_r\} \mbox{ and } x_l^{w(x_l)}x_i^{w(x_i)}x_j^{w(x_j)}x_r^{w(x_r)} \mid {\bf x}^{\bf a}.
\end{align*}
Since $x_u^{w(x_u)}x_l^{w(x_l)} \in I$ and by assumption $V^+$ are sinks, we have $w(x_u)=1$ or $w(x_l)=1$. Without loss of generality, suppose $w(x_u)=1$. This implies that $u \not \in \supp \mathbf{a}$, because ${\bf x}^{\bf a} \not\in I^{(3)}$. 
 \end{proof}
 \noindent 
 Now, we prove the first main result of this section. 
\begin{theorem}\label{thm1}
 Let $I=I(D)$ be an edge ideal of a weighted oriented graph $D$ such that $V^+$ are sink vertices. Then,
 \begin{equation*}
     \reg(I^{(k)}) \leq \reg(I^k) \mbox{ for } k=2,3.
 \end{equation*}
\end{theorem}
\begin{proof} Let $k=2$. Let ${\bf x}^{\bf a} \notin I^{(2)}$. Then, by Lemma $\ref{lm4.1}$, we have $\sqrt{I^{(2)}:{\bf x}^{\bf a}} = \sqrt{I^{2}:{\bf x}^{\bf a}}$. This implies that $\Delta_{\bf a}(I^{(2)})=\Delta_{\bf a}(I^{2})$. Therefore, by Lemma $\ref{lm2.1}$  we have  $\reg(I^{(2)}) \leq \reg(I^2)$, as required. Now, assume $k=3$. We show the required inequality by induction on $n=|V(D)|$. If $n=2$, then the required inequality is true because $I(D)^{(3)}=I(D)^3$. Assume $n\geq 3$. Let $(\mathbf{a},i)$ be an extremal exponent of $I^{(3)}$. Then, ${\bf x}^{\bf a} \notin I^{(3)}$. If $\Delta_{\mathbf{a}}(I^{(3)}) = \Delta_{\mathbf{a}}(I^3)$, then by Lemma \ref{lm2.1}, the required inequality follows. Assume $\Delta_{\mathbf{a}}(I^{(3)}) \neq \Delta_{\mathbf{a}}(I^3)$. By Lemma \ref{lm2.2} and Lemma \ref{lm4.2}, there exists a variable $x_j$ such that $x_j \in \sqrt{I^{(3)}:{\bf x}^{\bf a}}$ and $j \not \in \supp{\mathbf{a}}$. Let $U=[n]\setminus \{j\}$ and $J= I_{U}$. Then, by Lemma \ref{lm2.3}, we have  $\reg(I^{(3)})=\reg(J^{(3)})$. Also, by induction, we have $\reg(J^{(3)}) \leq \reg(J^3)$. Thus, we have $\reg(I^{(3)}) \leq \reg(J^3)$. But by Lemma \ref{lm2.4}, and the fact that $J^3$ is the restriction of $I^3$ to $U$, we have that $\reg(J^3) \leq  \reg(I^3)$. This implies that $\reg(I^{(3)}) \leq \reg(I^3)$, as required.
\end{proof}
\noindent
{\bf Notation:} Let  $\mathcal{F} $ be a family of all sets $\{x_i,x_j,x_r\}$ such that the induced subgraph on $\{x_i,x_j,x_r\}$ in $G$ is a triangle. Define $N[H] := \displaystyle \bigcup_{x \in V(H)} N[x]$, for any subgraph $H$ of $D$. 

\begin{proposition}\label{thm3.4}
Let $D$ be a weighted oriented graph such that $V^+$ are sink vertices and $I=I(D)$.  Then
\begin{align*}
 \reg(I^2) \leq &  \max \bigg \{\reg(I^{(2)}), \reg(I(D\setminus N[\supp T]))+
 \sum_{x \in N[\supp{T}]}w(x)-\vert N[\supp{T}]\rvert+ \\
 & \sum_{x \in \supp{T}} w(x) ~|~ \supp(T) \in \mathcal{F} \bigg \}.
\end{align*}
\end{proposition}
\begin{proof}
 By Lemma \ref{lm2.9}, we have that 
 $$I(G)^{(2)}=I(G)^2+(x_ix_jx_r \mid \{x_i,x_j,x_r\} \in \mathcal{F}).$$ Thus, by Lemma \ref{lm2.8} this implies that 
\begin{eqnarray*}
  I^{(2)} &=& \Phi(I(G)^{(2)}) \\
          &=& \Phi(I(G)^2)+(\Phi(x_ix_jx_r) \mid \{x_i,x_j,x_r\} \in \mathcal{F}) \\ 
        &=& I^2 +(x_i^{w(x_i)}x_j^{w(x_j)}x_r^{w(x_r)} \mid \{x_i,x_j,x_r\} \in \mathcal{F}).           
\end{eqnarray*}
Let $\{T_1,\ldots,T_t\}=\{(\Phi(x_ix_jx_r)) \mid \{x_i,x_j,x_r\} \in \mathcal{F} \}$ be the collection of all principal ideals such that each ideal is generated by monomial $x_i^{w(x_i)}x_j^{w(x_j)}x_r^{w(x_r)}$, for $\{x_i,x_j,x_r\} \in \mathcal{F}$. 
\noindent
First we show that $((I^2,T_1,\ldots,T_{i-1}):T_i)=(I^2:T_i)$ for each $i$. Note that $(I^2:T_i) \subseteq ((I^2,T_1,\ldots,T_{i-1}):T_i)$. Let $f \in ((I^2,T_1,\ldots,T_{i-1}):T_i)$. Then $fT_i \subseteq (I^2,T_1,\ldots,T_{i-1})$. If $fT_i \subseteq I^2$, then we are done. Suppose $fT_i \subseteq T_j$ for some $1 \leq j \leq i-1$. Then $fx_l^{w(x_l)}x_m^{w(x_m)}x_n^{w(x_k)}=gx_a^{w(x_a)}x_b^{w(x_b)}x_c^{w(x_c)}$, where $T_i= (x_l^{w(x_l)}x_m^{w(x_m)}x_k^{w(x_k)})$ and $T_j=(x_a^{w(x_a)}x_b^{w(x_b)}x_c^{w(x_c)})$ and $g \in R$. Suppose $\{x_l^{w(x_l)},x_m^{w(x_l)},x_k^{w(x_k)}\} \cap \{x_a^{w(x_a)},x_b^{w(x_b)},x_c^{w(x_c)}\}=\emptyset$. Then $x_a^{w(x_a)}x_b^{w(x_b)}x_c^{w(x_c)} \vert f$ and hence $fT_i \subseteq I^2$, as required. If $\{x_l^{w(x_l)},x_m^{w(x_m)},x_k^{w(x_k)}\} \cap \{x_a^{w(x_a)},x_b^{w(x_b)},x_c^{w(x_c)}\}=\{x_a^{w(x_a)}\} $, then $x_b^{w(x_b)}x_c^{w(x_c)} \vert f$ and hence $fT_i \subseteq I^2$ as required. If $\{x_l^{w(x_l)},x_m^{w(x_m)},x_k^{w(x_k)}\} \cap \{x_a^{w(x_a)},x_b^{w(x_b)},x_c^{w(x_c)}\} =\{x_a^{w(x_a)},x_b^{w(x_b)}\}$ and hence $fT_i \subseteq I^2$, as required. Now, consider the exact sequences:
\begin{align}
\begin{split}
 0 \rightarrow  (I^2:T_1)(-w_1) \rightarrow & I^2 \rightarrow (I^2+ T_1) \rightarrow 0  \\ 
\quad \; \vdots \hspace{2.5cm} &  \vdots \hspace{1.8cm}  \vdots\\
 0 \rightarrow  ((I^2,T_1.\ldots,T_{t-2}):T_{t-1}))(-w_{t-1}) \rightarrow & I^2+(T_1,\ldots,T_{t-2}) \rightarrow (I^2,T_1,\ldots,T_{t-1}) \rightarrow 0 \\
 0 \rightarrow  ((I^2,T_1.\ldots,T_{t-1}):T_t))(-w_t) \rightarrow & I^2+(T_1,\ldots,T_{t-1}) \rightarrow I^{(2)} \rightarrow 0,  \label{eq3}
 \end{split}
\end{align}
where $w_l= \sum_{x \in \supp{T_l} }w(x)$ for $1\leq l \leq t $.
Then by Lemma \ref{lm2.12}(1), we have
\begin{align*}
    \reg(I^2) \leq \max\{\reg(I^{(2)}), \reg((I^2:T_1))+w_1,\ldots,\reg((I^2,T_1,\ldots,T_{t-1})):T_t)+w_t\}.
\end{align*}
Therefore,
\begin{equation*}
    \reg(I^2) \leq \max \left \{\reg(I^{(2)}),\reg(I^2: x_i^{w(x_i)}x_j^{w(x_j)}x_r^{w(x_r)})+ w(x_i)+w(x_j)+w(x_r) \mid \{x_i, x_j, x_r\} \in \mathcal{F} \right \}. 
\end{equation*}

Now, we have 
\begin{equation*}
    (I^2:T_i)=I(D\setminus N[\supp T_i])+(x^{w(x)} \mid x \in N[\supp T_i]) \mbox{ for each } T_i.
\end{equation*}
Then, by \cite[Proposition 3.11.]{nv19}, for all $1\leq i \leq  t$, we have 
\begin{equation}\label{eq4}
    \reg((I^2:T_i))= \reg(I(D\setminus N[\supp T_i]))+\sum_{x \in N[\supp{T_i}]}w(x)-\vert N[\supp{T_i}]\rvert.
\end{equation}
 Thus, we have 
\begin{align*}
 \reg(I^2) \leq &  \max \bigg \{\reg(I^{(2)}), ~~\reg(I(D\setminus N[\supp T]))+
 \sum_{x \in N[\supp{T}]}w(x)-\vert N[\supp{T}]\rvert+ \\
 & \sum_{x \in \supp{T}} w(x) ~|~ \supp(T) \in \mathcal{F} \bigg \}.
\end{align*}
\end{proof}

 Below, we prove the next main result of this section. 
\begin{theorem}\label{thm3.5}
Let $D$ be a weighted oriented graph such that $V^+$ are sink vertices and $I=I(D)$.  Then
\begin{align*}
 \reg(I^2) \in &   \Bigg \{\reg(I^{(2)}), \max\bigg\{\reg(I(D\setminus N[\supp T]))+
 \sum_{x \in N[\supp{T}]}w(x)-\vert N[\supp{T}]\rvert+ \\
 & \sum_{x \in \supp{T}} w(x) ~|~ \supp(T) \in \mathcal{F}\bigg\} \Bigg\}.
\end{align*}
\end{theorem}

\begin{proof}
If $\reg(I^2)=\reg(I^{(2)})$, then we are done. Now, assume  $\reg(I^2) \neq \reg(I^{(2)})$.  Then by Theorem \ref{thm1}, we have $\reg(I^{(2)}) < \reg(I^2)$. Assume the notations as in the proof of Proposition \ref{thm3.4}. Without loss of generality, assume that 
\begin{equation} \label{eq5}
 \max\{ \reg((I^2:T_1))+w_1,\ldots,\reg((I^2:T_{t-1}))+w_{t-1},\reg((I^2:T_{t}))+w_{t} \}=\reg((I^2:T_t))+w_t.   
\end{equation}
Note that from proof of Proposition \ref{thm3.4}, we have $((I^2,T_1,\ldots,T_{i-1}):T_t)=(I^2:T_i)$ for each $i$. Now, consider the exact sequences as in \eqref{eq3} and by Lemma \ref{lm2.12}(2), we have
\begin{align*}
 \reg(((I^2,T_1,\ldots,T_{t-1}):T_t))+w_t \leq \max\{ \reg((I^2,T_1,\dots,T_{t-1})), \reg(I^{(2)})+1\}.
\end{align*}
Now, by applying Lemma \ref{lm2.12}(3) recursively to exact sequence \eqref{eq3}, we have
\begin{align*}
    \reg((I^2,T_1,\dots,T_{t-1})) \leq & \max \bigg \{\reg(I^2), 
      \reg((I^2:T_1))+w_1-1,\ldots,  \\
      & \reg((I^2:T_{t-1}))+w_{t-1}-1 \bigg \} 
\end{align*}
which results in 
\begin{align*}
     \reg((I^2:T_t))+w_t = & \reg(((I^2,T_1,\ldots,T_{t-1}):T_t))+w_t   \\
      \leq & \max \bigg \{\reg(I^{(2)})+1, \reg(I^2), 
      \reg((I^2:T_1))+w_1-1,  \\
      & \ldots,\reg((I^2:T_{t-1}))+w_{t-1}-1 \bigg \} \\
      \leq & \max \{\reg(I^{(2)})+1, \reg(I^2) \} (\mbox{ by equation } \eqref{eq5}) \\
      \leq & \reg(I^2). 
\end{align*}
Therefore, from the equation \eqref{eq5}, we get that  
\begin{align*}
    \reg((I^2:T_i))+w_i \leq \reg(I^2), \mbox{ for all } i.   
\end{align*}

Thus, by equation \eqref{eq4}, for each $T_i$, we have
\begin{equation*}
\reg(I^2) \geq \reg(I(D\setminus N[\supp T_i]))+\sum_{x \in N[\supp{T_i}]}w(x)-\vert N[\supp{T_i}]\rvert+\sum_{x \in \supp{T_i}} w(x).
\end{equation*}
Also, by Proposition \ref{thm3.4} and $\reg(I^{(2)}) < \reg(I^2)$, we have 
\begin{align*}
 \reg(I^2) \leq & \max \bigg \{ \reg(I(D\setminus N[\supp T]))+
 \sum_{x \in N[\supp{T}]}w(x)-\vert N[\supp{T}]\rvert+ \\
 & \sum_{x \in \supp{T}} w(x) ~|~ \supp(T) \in \mathcal{F} \bigg\}.
\end{align*}
Thus, we have
\begin{align*}
 \reg(I^2) = & \max \bigg \{ \reg(I(D\setminus N[\supp T]))+
 \sum_{x \in N[\supp{T}]}w(x)-\vert N[\supp{T}]\rvert+ \\
 & \sum_{x \in \supp{T}} w(x) ~|~ \supp(T) \in \mathcal{F} \bigg \},
\end{align*}
as required.
\end{proof}

\begin{corollary}\label{cor3.6}
Let $D$ be a weighted oriented gap-free graph with $V^+$ are sink vertices and $I=I(D)$. Then 
\begin{equation*}
 \reg(I^2)\leq  \displaystyle  \max \bigg\{\reg(I^{(2)}), \sum_{x \in N[\supp{T}]}w(x)-\vert N[\supp{T}]\rvert+1+\sum_{x \in \supp{T}} w(x) ~|~ \supp(T) \in \mathcal{F} \bigg\}.
\end{equation*}
\end{corollary}
\begin{proof}
 Let $G$ be a gap-free graph. If possible, there exists an edge $\{x,y\}$ of $ G \setminus N[\supp T]$, then the edge $\{x,y\}$ and any edge of the triangle $T$ forms a gap which is a contradiction. Thus, $D\setminus N[\supp T]$ is an empty graph or it has only isolated vertices. Thus, by Proposition \ref{thm3.4}, we have
\begin{equation*}
 \reg(I^2)\leq  \displaystyle  \max\bigg \{\reg(I^{(2)}), \sum_{x \in N[\supp{T}]}w(x)-\vert N[\supp{T}]\rvert+1+\sum_{x \in \supp{T}} w(x) ~|~ \supp(T) \in \mathcal{F} \bigg\}.
\end{equation*}
\end{proof}
\begin{corollary}\label{cor3.8}
Let $D$ be a weighted oriented gap-free graph with $V^+$ are sink vertices. Let $I=I(D)$. Then 
\begin{equation*}
 \reg(I^2) \in  \displaystyle \Bigg\{\reg(I^{(2)}), \max \bigg\{ \sum_{x \in N[\supp{T}]}w(x)-\vert N[\supp{T}]\rvert+1+\sum_{x \in \supp{T}} w(x) ~|~ \supp(T) \in \mathcal{F} \bigg\}\Bigg\}.
\end{equation*}
\end{corollary}
\begin{proof}
Note that $D\setminus N[\supp T]$ is an empty graph or it has only isolated vertices as $D$ is a weighted oriented gap-free graph. Therefore, by Theorem \ref{thm3.5}, we have 
\begin{equation*}
 \reg(I^2) \in  \displaystyle \Bigg\{\reg(I^{(2)}), \max \bigg\{ \sum_{x \in N[\supp{T}]}w(x)-\vert N[\supp{T}]\rvert+1+\sum_{x \in \supp{T}} w(x) ~|~ \supp(T) \in \mathcal{F} \bigg\}\Bigg\},
\end{equation*}
 as required. 
\end{proof}
\begin{corollary}\label{cor3.7}
Let $D$ be a weighted oriented gap-free graph with $V^+$ are sink vertices and $I=I(D)$. For any set $\{x_i,x_j,x_r\} \in \mathcal{F}$, suppose $\vert N[\{x_i,x_j,x_r\}] \cap V^+ \rvert \leq 1$. Then $\reg(I^2)=\reg(I^{(2)})$.
\end{corollary}
\begin{proof}
 If $N[\{x_i,x_j,x_r\}] \cap V^+ \rvert \leq 1$ then  
 \begin{equation*}
  \displaystyle \sum_{x \in N[\supp{T}]}w(x)-\vert N[\supp{T}]\rvert+1+\sum_{x \in \supp{T}} w(x) \leq 2(w+1),   
 \end{equation*}
where $w=\max\{w(x) \mid x \in V(D)\}$. Since $(x_ix_j^{w(x_j)})^2$ is a minimal generator of $I^{(2)}$ for any edge $(x_i,x_j)\in E(D)$,  we get that $2(w+1)$ is less than or equal to the maximal degree of the minimal generators of $I^{(2)}$. Therefore, by Corollary \ref{cor3.6}, we have $\reg(I^2) \leq \max\{\reg(I^{(2)}), 2(w+1)\} \leq \reg(I^{(2)})$. Thus, by Theorem \ref{thm1}, we have $\reg(I^2)=\reg(I^{(2)})$, as required.
\end{proof}
\noindent
Now, we give some examples in which we see that the reverse inequality $\reg(I^k) \leq \reg(I^{(k)})$ will not hold in general and also depend on weights for $k=2,3$.
\begin{example}\label{notequalsink}
Let $D$ be a Cohen-Macaulay weighted oriented graph with all leaves are sinks ( $V^+$ are sinks ) as in Figure \ref{fig2}.
\begin{figure}[!h]
    \centering
     \includegraphics[width=0.35 \textwidth]{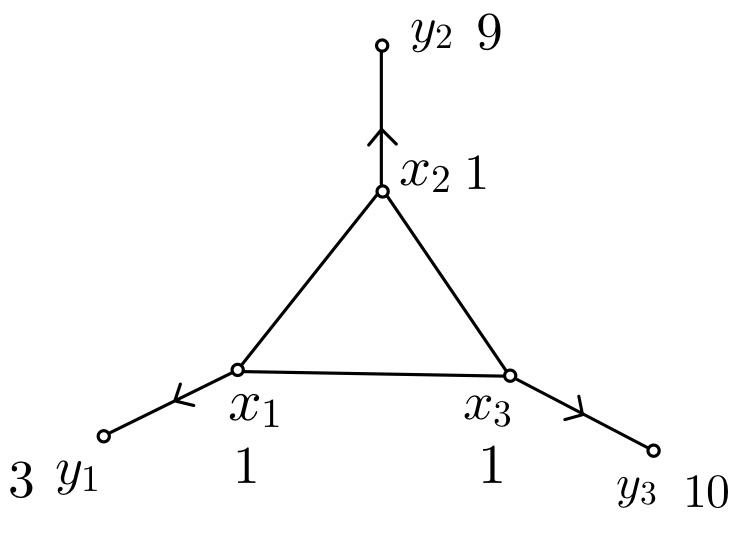} 
    \caption{Cohen-Macaulay weighted oriented graph}
    \label{fig2}
\end{figure}
Then, $ I(D)=(x_1x_2, x_2x_3,x_3x_1,x_1y_1^{3}, x_2y_2^{9}, x_3y_3^{10}).$ \\
Using Macaulay2 \cite{gs}, we have
\begin{align*}
 22=\reg(I(D)^{(2)}) <& \reg(I(D)^{2})=23 \\
 33=\reg(I(D)^{(3)}) <& \reg(I(D)^{3})=34.   
\end{align*}
\end{example} 
\begin{example}
Let $D$ be a Cohen-Macaulay weighted oriented graph with all leaves are sinks ( $V^+$ are sinks ) as in Figure \ref{fig3}.
\begin{figure}[!h]
    \centering
     \includegraphics[width=0.35 \textwidth]{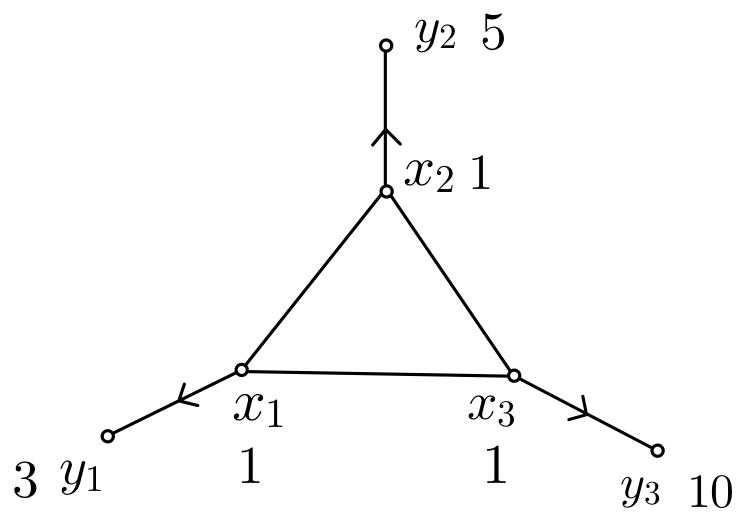} 
    \caption{Cohen-Macaulay weighted oriented graph}
    \label{fig3}
\end{figure}
Then, $I(D)=(x_1x_2, x_2x_3,x_3x_1,x_1y_1^{3}, x_2y_2^{5}, x_3y_3^{10}).$ \\
Using Macaulay2 \cite{gs}, we have 
\begin{align*}
 \reg(I(D)^{(2)}) = 22 &= \reg(I(D)^{2}) \\
 \reg(I(D)^{(3)}) = 33 &= \reg(I(D)^{3}).   
\end{align*}
\end{example}  

\noindent 
The below example shows that the inequality $\reg(I(D)^{(k)}) \leq \reg(I(D)^k)$ can be strict even in case $V^+$ are not sinks.
\begin{example}\label{notequalnotsink}
	Let $D$ be a Cohen-Macaulay weighted oriented forest as in Figure \ref{fig4}.
\begin{figure}[!h]
    \centering
     \includegraphics[width=0.25 \textwidth]{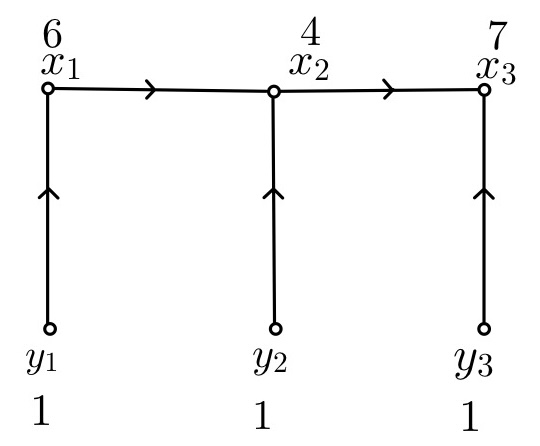} 
    \caption{Cohen-Macaulay weighted oriented forest}
    \label{fig4}
\end{figure}
Then, $ I(D)=(x_1x_2^4,x_2x_3^7,x_1^6y_1,y_2x_2^4,y_3x_3^7).$ 
Using Macaulay2 \cite{gs}, we have 
\begin{align*}
   \reg(I(D)^{(2)})=23 <& \reg(I(D)^{2})=24 \\
   \reg(I(D)^{(3)})=30 <&\reg(I(D)^{3})=32. 
\end{align*}
\end{example}


\section{Regularity of symbolic and ordinary powers of complete graph }\label{sec5}

In this section, we give explicit formulas for the regularity of symbolic and ordinary powers of edge ideals of a weighted oriented complete graph $K_n$. For simplicity, we denote the ideals $I_{C}^{K_n}$ by $I(C)$ for any $C\in \mathcal{C}_s$, where $\mathcal{C}_s$ is the set of all strong vertex cover of $K_{n}$. Note that $C_i=N_{K_n}(x_i)$ and $C_j=N_{K_n}(x_j)$. Also, $L_1(C_i)=N^-_{K_n}(x_i)$ and $L_2(C_i)=N^+_{K_n}(x_i)$. \\ 

\noindent 
Now, we state a lemma whose proof follows from \cite[Corollary 3.4]{b23}. 

\begin{lemma}\label{lem-Cohen-Macaulay}
  Let $K_{n}$ be a weighted oriented complete graph such that $N^+_{K_n}(V^+) \neq V(K_n)$.  Then, $I(K_n)^{(k)}$ is Cohen-Macaulay, for any $k\geq 1$.  
\end{lemma}

Now, we give a remark which will be used to prove the main results of this section. 
\begin{remark} \label{rmk1}
$({\bf 1})$ If $\Delta=\{\emptyset\}$ is an empty complex then for any $F\in \Delta$, we have $\mathrm{lk}_{\Delta} F$ is an empty complex. \\
$({\bf 2})$  If $\Delta$ is a single point, that is, $\Delta=\Big\{\emptyset,\{i\}\Big\}$, for some $i\in [n]$, then for any $F\in \Delta$, we have $\mathrm{lk}_{\Delta} F$ is an empty complex or a single point.  \\
$({\bf 3})$  If $\Delta=\Big\{\emptyset, \{i_1\},\{i_2\},\ldots, \{i_r\}\Big\}$ is a set of isolated points and $r\geq 2$, then for any $F\in \Delta$, we have $\mathrm{lk}_{\Delta} F$ is an empty complex or $\Delta$.  \\ 
$({\bf 4})$ For a complete graph $K_n$, $\Delta(I(K_n))=\Big\{ \emptyset, \{1\}, \{2\}, \ldots,\{n-1\},\{n\}\Big\}$. Let $k\geq 1$. For $\mathbf{a} = (a_1, \ldots, a_n) \in \NN^n$, the degree complexes $\Delta_{\mathbf{a}}(I(K_n)^{(k)})$ or $\Delta_{\mathbf{a}}(I(K_n)^{k})$ are subcomplexes of $\Delta(I(K_n))$. Then $\Delta_{\mathbf{a}}(I(K_n)^{(k)})$ or $\Delta_{\mathbf{a}}(I(K_n)^{k})$ is a single point or a set of at least two isolated points or the empty complex $\{\emptyset \}$. \\
$({\bf 5})$ Let \(\mathbf{a} \in \mathbb{N}^n\) be an exponent. Then \(\{i\}\) is a face of \(\Delta_{\mathbf{a}}(I(K_n)^{(k)})\) (or \(\Delta_{\mathbf{a}}(I(K_n)^{k})\) if and only if \(x^{\mathbf{a}} \not\in I(C_i)^k\), which is equivalent to
 $\displaystyle \sum_{x_l \in N^-_{K_n}(x_i)} a_l + \sum_{x_l \in N^+_{K_n}(x_i)} \left\lfloor \frac{a_l}{w(x_l)} \right\rfloor < k$, where $C_i=V(K_n)\setminus \{x_i\}$. 
\end{remark}
\noindent 
{\bf Convention:} We set  $\max\{w(x) \mid x \in A\} = 1$ if $A$ is an empty subset of $V(K_n)$.
\noindent
\begin{lemma} \label{lm-max1}
For each \((x_i, x_j) \in E(K_n)\) and each positive integer \(k\), let \(A(i, j, k)\) be the set of exponents \(\mathbf{a} = (a_1, \ldots, a_n) \in \mathbb{N}^n\) such that
$$\sum_{x_l \in N^-_{K_n}(x_i)} a_l + \sum_{x_l \in N^+_{K_n}(x_i)} \left\lfloor \frac{a_l}{w(x_l)} \right\rfloor < k  \text{ and } \sum_{x_l \in N^-_{K_n}(x_j)} a_l + \sum_{x_l \in N^+_{K_n}(x_j)} \left\lfloor \frac{a_l}{w(x_l)} \right\rfloor < k.$$
Then, 
 \begin{align*}
      r_{ij}(k) =& \max\{|\mathbf{a}| \mid \mathbf{a} \in A(i, j, k)\} \\
     = & \max\{w(x_j)+1,w(x) \mid x \in  N^+_{K_n}(x_i)\cap N^+_{K_n}(x_j)\}(k-1)+\sum_{x \in N^+_{K_n}(x_i)
      \cap N^+_{K_n}[x_j]}(w(x)-1).
 \end{align*}
 \end{lemma}
\begin{proof}
For each \((x_i, x_j) \in E(K_n)\), we define 
\[
B_1 = N^-(x_i) \cap N^-(x_j), B_2 = N^-(x_i) \cap N^+(x_j), B_3 = N^+(x_i) \cap N^-(x_j), B_4 = N^+(x_i) \cap N^+(x_j).
\]
Then, we have \(B_1, B_2, B_3, B_4\) are disjoint and 
\(B_1 \cup B_2 \cup B_3 \cup B_4 =  V(K_n) \setminus \{x_i, x_j\}.\)
The given inequalities can be written as follows: 
\begin{align*}
 &\sum_{x_t \in B_1} a_t +\sum_{x_t \in B_2} a_t+ \sum_{x_t \in B_3} \left\lfloor \frac{a_t}{w(x_t)} \right\rfloor +\sum_{x_t \in B_4} \left\lfloor \frac{a_t}{w(x_t)} \right\rfloor +\left\lfloor \frac{a_j}{w(x_j)} \right\rfloor < k   \\  &\text{ and } \\
&\sum_{x_t \in B_1} a_t +\sum_{x_t \in B_2} \left\lfloor \frac{a_t}{w(x_t)} \right\rfloor+ \sum_{x_t \in B_3} a_t +\sum_{x_t \in B_4} \left\lfloor \frac{a_t}{w(x_t)} \right\rfloor+a_i < k.
\end{align*}
Now, assume that \(\mathbf{a} \in A(i, j, k)\) is an exponent satisfying the above two inequalities. Now, we perform the following reduction: If \(a_l \neq 0\) for some \(l \in B_1\), then consider \(\mathbf{a}' \in \mathbb{N}^n\) defined by setting 
\(
a'_l = 0, \quad a'_j = a_j + a_l, \quad a'_i = a_i + a_l.\)
Then, we have \(\mathbf{a}' \in A(i, j, k)\) and \(|\mathbf{a}'| \geq |\mathbf{a}|\). Similarly, the same reduction applies for the case \(a_l \neq 0\) for some \(l \in B_2\) or \(B_3\). Hence, we may assume that the extremal solution for \(r_{ij}(k)\) is such that \(a_t = 0\) for all \(t \in B_1 \cup B_2 \cup B_3\). Then, we have following reduced inequalities
\begin{align*}
 \sum_{x_t \in B_4} \left\lfloor \frac{a_t}{w(x_t)} \right\rfloor +\left\lfloor \frac{a_j}{w(x_j)} \right\rfloor < k    \text{ and } \sum_{x_t \in B_4} \left\lfloor \frac{a_t}{w(x_t)} \right\rfloor+a_i < k.
\end{align*}
Now, suppose \( w(x_u) = \max \{w(x) \mid x \in B_4\} \) and consider two cases: \\
{\bf Case 1:} Suppose \(w(x_j) < w(x_u)\). Then, for extremal solution, we have $a_u=w(x_u)k-1$, $a_j=w(x_j)-1$, $a_i=0$, and $a_t=w(x_t)-1$ for all $x_t \in B_4\setminus \{x_u\}$. Thus, we have
\begin{align*}
  \max\{|\mathbf{a}| \mid \mathbf{a} \in A(i, j, k)\}=&w(x_u)k-1+(w(x_j)-1)+\sum_{x_t \in B_4\setminus \{x_u\}}(w(x_t)-1) \\
  =& w(x_u)(k-1)+(w(x_u)-1)+(w(x_j)-1)+\sum_{x_t \in B_4\setminus \{x_u\}}(w(x_t)-1) \\
  =& w(x_u)(k-1)+\sum_{x \in N^+_{K_n}(x_i) \cap  N^+_{K_n}[x_j]} (w(x)-1).
\end{align*}
\noindent
{\bf Case 2:} Suppose $w(x_j) \geq w(x_u)$. Then, for extremal solution, we have $a_j=w(x_j)k-1$, $a_i=k-1$, $a_t=w(x_t)-1$ for all $x_t \in B_4$. Thus, we have
\begin{align*}
  \max\{|\mathbf{a}| \mid \mathbf{a} \in A(i, j, k)\}=&w(x_j)k-1+(k-1)+\sum_{x_t \in B_4}(w(x_t)-1) \\
  =& (w(x_j)+1)(k-1)+(w(x_j)+1)-2+\sum_{x_t \in B_4}(w(x_t)-1) \\
  =& (w(x_j)+1)(k-1)+\sum_{x \in N^+_{K_n}(x_i) \cap  N^+_{K_n}[x_j]} (w(x)-1).
\end{align*}
Thus, by combining {\bf Case 1:} and {\bf Case 2:}, we have 
\begin{align*}
r_{ij}(k)=\max \{w(x_j) + 1, w(x) \mid x \in N^+_{K_n}(x_i) \cap N^+_{K_n}(x_j)\}(k-1) + \sum_{x \in N^+_{K_n}(x_i) \cap N^+_{K_n}[x_j]} (w(x) - 1).
\end{align*}
\end{proof} 
The proof of below lemma is simple and can be omitted.
\begin{lemma}\label{lm-max2}
For each \(i \in [n]\) and each positive integer \(k\), let \(P(i, k)\) be the set of exponents \(\mathbf{a} = (a_1, \ldots, a_n)\) such that
$$\sum_{x_l \in  N^-_{K_n}(x_i)} a_l+\sum_{x_l \in  N^+_{K_n}(x_i)} \left\lfloor \frac{a_l}{w(x_l)} \right\rfloor < k. $$ 
Then,
\begin{align*}
r_{i}(k)=\max\{|\mathbf{a}| \mid \mathbf{a} \in P(i,k)\}=\max\{w(x) \mid x \in N^+_{K_n}(x_i)\}(k-1)+\sum_{x \in N^+_{K_n}(x_i)}(w(x) -1).
\end{align*}

\end{lemma}
\noindent 
 
Now, we prove the first main result of this section.

\begin{theorem}\label{thmcompletesymbolic1}
 Let $K_n$ be a weighted oriented complete graph such that \( N^+_{K_n}(V^+) \neq V(K_n) \). Then, for any $k\geq 1$,  
\begin{align*}
\reg(I(K_n)^{(k)}) = \max \bigg\{\underset{ (x_i,x_j) \in E(K_n)}{\max}\{ r_{ij}(k)\} + 2, 
\underset{x_i \in V(K_n)}{\max} \{r_i(k)\} + 1 \bigg\}, \text{ where }
\end{align*}
$$r_{ij}(k)= \max \{w(x_j) + 1, w(x) \mid x \in N^+_{K_n}(x_i) \cap N^+_{K_n}(x_j)\}(k-1) + \sum_{x \in N^+_{K_n}(x_i) \cap N^+_{K_n}[x_j]} (w(x) - 1)$$ 
and $r_{i}(k)=\max\{w(x) \mid x \in N^+_{K_n}(x_i)\}(k-1)+\sum_{x \in N^+_{K_n}(x_i)}(w(x) -1).$
\end{theorem}
\begin{proof}
We prove the required formula using Lemma \ref{hom1}.
Let \((\mathbf{a}, i)\) be an extremal exponent of \(I(K_n)^{(k)}\). Since \(\Delta_{\mathbf{a}}(I(K_n)^{(k)})\) is a subsimplicial complex of \(\Delta(I(K_n)) = \{\{1\}, \ldots, \{n\}\}\), we deduce that \(i = 1\) or \(i = 0\). Now, assume that \(i = 1\). Then 
\(
\widetilde{H}_0(\Delta_{\mathbf{a}}(I(K_n)^{(k)});\K) \neq 0.
\)
This implies that there must exist \(i, j\) such that $\{i\}, \{j\} \in \Delta_{\mathbf{a}}(I(K_n)^{(k)})$. Hence, by Remark \ref{rmk1}(5) and Lemma \ref{lm-max1}, we get the first term.

Now, assume that \(i = 0\). By Lemma \ref{lem-Cohen-Macaulay},  \(\mathfrak{m}\) is not an associated prime of \(I(K_n)^{(k)}\), we deduce that \(\Delta_{\mathbf{a}}(I(K_n)^{(k)})\) is not the empty complex. This implies that \(F = \{l\}\) for some \(l \in [n]\). In particular, \(\{l\}\) is a face of \(\Delta_{\mathbf{a}}(I(K_n)^{(k)})\). Hence, by Remark \ref{rmk1}(5) and Lemma \ref{lm-max2}, we get the second term, and the conclusion follows, by Lemma \ref{hom1}.
\end{proof}

Note that by \cite[Corollary 3.4]{b23}, the condition $N^+_{K_n}(V^+) \neq V(K_n)$ is equivalent to  $I(K_n)^{(k)} \neq I(K_n)^k$ for all $k\geq 2$. Also, we can not drop this condition in Theorem \ref{thmcompletesymbolic1}. See the below example. 
\begin{example}
Consider a weighted oriented complete graph with edge ideal $I(K_4)=(x_1x_2^3,x_1x_4^6,x_3x_1^9,x_3x_2^3,x_4x_3^5,x_2x_4^6)$. One can see that $N^{+}_{K_4}(V^+)=V(K_4)$. 
By Macaulay2 \cite{gs}, we have $$\reg(I(K_4)^{(2)})=\reg(I(K_4)^2)=30.$$
Note that
\begin{align*}
& r_{12}(k)=6(k-1)+7, r_{14}(k)=7(k-1)+5, r_{31}(k)=10(k-1)+10 , \\
& r_{32}(k)=4(k-1)+2, r_{43}(k)=6(k-1)+4, r_{24}(k)=7(k-1)+5. 
\end{align*}
Also, we have
\begin{align*}
      r_{1}(k)=6(k-1)+7, ~~~r_2{(k)}=6(k-1)+5, ~~~r_3(k)=9(k-1)+10, ~~~r_4(k)=5(k-1)+4.
\end{align*}
Thus, the formula in Theorem \ref{thmcompletesymbolic1}, for any $k\geq 2$, 
$$\max \bigg ( \underset{ (x_i,x_j) \in E(K_4)}{\max}\{ r_{ij}(k)\} + 2, 
\underset{x_i \in V(K_4)}{\max} \{r_i(k)\} + 1 \bigg )=10(k-1)+12.$$ 
In particular, for $k=2$, we have
$$\max \bigg ( \underset{ (x_i,x_j) \in E(K_4)}{\max}\{ r_{ij}(2)\} + 2, 
\underset{x_i \in V(K_4)}{\max} \{r_i(2)\} + 1 \bigg )=22 \neq \reg(I(K_4)^{(2)}).$$ 
\end{example}
\noindent 
Below we give a formula for the regularity of $I(K_n)^k$ and comparison with that of symbolic powers. 

\begin{theorem}\label{thmcompletesymbolic2}
 Let $K_{n}$ be a weighted oriented complete graph. Then for any $k\geq 1$, 
{\tiny 
 \begin{align*}
 \reg(I(K_n)^{k}) = \max \bigg\{\underset{ (x_i,x_j) \in E(K_n)}{\max}\{ r_{ij}(k)\} + 2, 
 \underset{x_i \in V(K_n)}{\max}\{ r_i(k) \} + 1, \underset{(\alpha_1,\ldots,\alpha_n) \in \mathcal{M}_{\m}}{\max} \Big\{\sum_{i=1}^n (\alpha_i-1) \Big\}+1   \bigg\}.
\end{align*}
}
Moreover, $\reg(I(K_n)^{(k)}) \leq \reg(I(K_n)^k), \text{ for all } k \geq 1.$
\end{theorem}
\begin{proof}
We prove the required formula using Lemma \ref{hom1}.
Let \((\mathbf{a}, i)\) be an extremal exponent of \(I(K_n)^{k}\). Since \(\Delta_{\mathbf{a}}(I(K_n)^{k})\) is a subsimplicial complex of \(\Delta(I(K_n)) = \{\{1\}, \ldots, \{n\}\}\), we deduce that \(i = 1\) or \(i = 0\). Now, assume that \(i = 1\). Then 
\(
\widetilde{H}_0(\Delta_{\mathbf{a}}(I(K_n)^{k});\K) \neq 0.
\)
This implies that there must exist \(i, j\) such that $\{i\}, \{j\} \in \Delta_{\mathbf{a}}(I(K_n)^{k})$. Hence, by Remark \ref{rmk1}(5) and Lemma \ref{lm-max1}, we get the first term.

Now, assume that \(i = 0\). If \(\Delta_{\mathbf{a}}(I(K_n)^{k})\) is not the empty complex then \(F = \{l\}\) for some \(l \in [n]\). In particular, \(\{l\}\) is a face of \(\Delta_{\mathbf{a}}(I(K_n)^{k})\). Hence, by Remark \ref{rmk1}(5) and Lemma \ref{lm-max2}, we get the second term.
Now, suppose $\Delta_{\mathbf{a}}(I(K_n)^k)=\{ \emptyset\}$, the empty complex. Then their corresponding ideals gives that
 \begin{align*}
(x_1,\ldots,x_n)= I_{\Delta_{\mathbf{a}}(I(K_n)^{k})} =\bigcap_{(\alpha_1,\ldots,\alpha_n)\in \mathcal{M}_{\m}} \sqrt{(x_1^{\alpha_1}, \ldots, x_n^{\alpha_n}): {\bf x}^{\bf a}}. 
\end{align*}
This implies that there exists $(\alpha_1,\ldots,\alpha_n) \in \mathcal{M}_{\m}$, ${\bf x}^{\bf a} \not \in (x_1^{\alpha_1}, \ldots, x_n^{\alpha_n})$  which is equivalent to the inequality
$  \left\lfloor \frac{a_1}{\alpha_1} \right\rfloor+ \cdots +\left\lfloor \frac{a_n}{\alpha_n} \right\rfloor < 1. $ Therefore, we get the third term and hence we get the required formula, by Lemma \ref{hom1}.
\par If $I(K_n)^{(k)} \neq I(K_n)^{k}$, then by Theorem \ref{thmcompletesymbolic1} and using the formula for $\reg(I(K_n)^k)$, we have
$\reg(I(K_n)^{(k)}) \leq \reg(I(K_n)^k)$. This finishes the proof. 
\end{proof}

\noindent 
As a consequence, we show that $\reg(I(K_n)^{(k)})$ is eventually a linear function of $k$. 
\begin{corollary}\label{eventually-linear}
Let $K_n$ be a weighted oriented complete graph. Then $\reg(I(K_n)^{(k)})$ is eventually a linear function of $k$. Moreover, if \( N^+_{K_n}(V^+) \neq V(K_n) \), then 
$$\reg(I(K_n)^{(k)}) = (\omega+1)(k-1)+\underset{ (x_i,x_{\ell}) \in E(K_n)}{\max}  \Bigg\{\sum_{x \in N^+(x_i) \cap N^+[x_{\ell}]} (w(x) - 1) \Bigg\} + 2,$$
for $k\gg 0$, where $\omega=\max\{w(x) \mid x \in V(K_n)\}=w(x_{\ell})$ and maximum is taken over all $x_i, x_{\ell}\in V(K_n)$ with $w(x_{\ell})=\omega$. 
\end{corollary}
\begin{proof}
If \( N^+_{K_n}(V^+) = V(K_n) \) then by \cite[Corollary 3.4]{b23}, we have $I(K_n)^{(k)}=I(K_n)^{k}$ for all $k \geq 1$. Thus, the corollary follows from \cite[Theorem 1.1]{CHT99}. Assume \( N^+_{K_n}(V^+) \neq V(K_n) \). Let $\omega=\max\{w(x) \mid x \in V(K_n)\}$. Then there exists an edge $(x_{s},x_{\ell}) \in E(K_n)$ with $w(x_{\ell})=\omega$. Then, we have 
$$\underset{ (x_i,x_j) \in E(K_n)}{\max} \Big\{ w(x_j) + 1, w(x) \mid x \in N^+_{K_n}(x_i) \cap N^+_{K_n}(x_j) \Big\}=\omega+1 \geq \underset{ x_i \in V(K_n)}{\max} \Bigg\{ \underset{x \in N^+(x_i)}{\max}\{ w(x)\} \Bigg\}.$$
Then by Theorem \ref{thmcompletesymbolic1}, we get that for $k\gg 0$, 
$$\reg(I(K_n)^{(k)}) = (\omega+1)(k-1)+\underset{ (x_i,x_{\ell}) \in E(K_n)}{\max}  \Bigg\{\sum_{x \in N^+(x_i) \cap N^+[x_{\ell}]} (w(x) - 1) \Bigg\} + 2.$$
Thus, $\reg(I(K_n)^{(k)})$ is eventually a linear function of $k$. 
\end{proof}

\begin{example}
Let $I(K_4)=(x_1x_2^9,x_1x_3^9,x_1x_4^4,x_2x_4^4,x_4x_3^9,x_3x_2^9) \subset R=\mathbb{Q}[x_1,x_2,x_3,x_4]$ be an edge ideal of weighted oriented complete graph $K_4$. One can check that $N^+_{K_4}(V^+)\neq V(K_4)$. By Macaulay2 \cite{gs}, $\reg(I(K_4)^{(2)})=29$. 
Note that for all $k\geq 2$, 
\begin{align*}
 r_{12}(k)= 10(k-1)+11, ~~~ r_{13}(k)= 10(k-1)+16, & ~~~r_{14}(k)= 9(k-1)+11,  \\
 r_{24}(k)= 5(k-1)+3,~~~ r_{43}(k)= 10(k-1)+8, & ~~~r_{32}(k)= 10(k-1)+8 . 
\end{align*}
Also, we have
\begin{align*}
    r_{1}(k)=9(k-1)+19, ~~~r_2{(k)}=4(k-1)+3, ~~~r_3(k)=9(k-1)+8, ~~~r_4(k)=9(k-1)+8.
\end{align*}
Thus, we have
$$\reg(I(K_4)^{(2)}) = \max \bigg ( \underset{ (x_i,x_j) \in E(K_4)}{\max}\{ r_{ij}(2)\} + 2, 
\underset{x_i \in V(K_4)}{\max} \{r_i(2)\} + 1 \bigg )=29.$$
Furthermore, $\reg(I(K_4)^{(k)}) = 10(k-1)+18$, for all $k\geq 4$. 
\end{example}

\begin{remark}
For any disjoint union of weighted oriented complete graphs \\
$D=K_{n_1} \cup \ldots \cup K_{n_l}$, one can get explicit formulas for $\reg(I(D)^{(k)})$ and $\reg(I(D)^{k})$ by using \cite[Corollary 4.6]{hjkn23}, Theorem \ref{thmcompletesymbolic1} and \cite[Proposition 5.1]{nv19}, Theorem \ref{thmcompletesymbolic2} respectively.     
\end{remark}
\noindent 
After checking several examples with Macaulay2 \cite{gs}, we raise the following question:
\begin{question}
Let $D$ be a weighted oriented graph. Is $\reg(I(D)^{(k)}) \leq \reg(I(D)^k)$ for all $k\geq 1$?   
\end{question}
\noindent

\noindent 
{\bf Acknowledgement:} The authors would like to thank the anonymous referee for careful reading of the manuscript and suggestions for improving the paper to the present version. Manohar Kumar is thankful to the Government of India for supporting him in this work through the Prime Minister's Research Fellowship.


\begin{thebibliography}{31}

\bibitem{b23} A. Banerjee, B. Chakraborty, K. K. Das, M. Mandal, and S. Selvaraja, {\em Equality of ordinary and symbolic powers of edge ideals of weighted oriented graphs}, Comm. Algebra {\bf 51} (2023), no.~4, 1575-1580.

\bibitem{CHT99} S. D. Cutkosky, J. Herzog\ and\ N. V. Trung, {\em Asymptotic behaviour of the Castelnuovo-Mumford regularity}, Compositio Math. {\bf 118} (1999), no.~3, 243--261.	

\bibitem{chhktt19} G. Caviglia, H. T.  Hà, J. Herzog, M. Kummini, N. Terai and N. Trung, {\em Depth and regularity modulo a principal ideal},  J. Algebraic Combin. \textbf{49} (2019), no.~1, 1-20.

\bibitem{ctv93} M. Catalisano, N. Trung and G. Valla, {\em A sharp bound for the regularity index of fat points in general position},  Proc. Amer. Math. Soc. \textbf{118} (1993), no.~3, 717-724. 

\bibitem{czw22} Y. Cui, G. Zhu and X. Wei, {\em The edge ideals of the join of some vertex weighted oriented graphs}, to appear in J. Algebra. Appl. \url{https://doi.org/10.1142/S0219498825502962} 

\bibitem{dhnt21} L. Dung, T. Hien, H. Nguyen and  T. Trung, {\em Regularity and Koszul property of symbolic powers of monomial ideals}, Math. Z. \textbf{298} (2021), no.~3-4, 1487-1522.

\bibitem{gmv21} G. Gonzalo, J. Martinez-Bernal, and R. H. Villarreal, {\em Normally torsion-free edge ideals of weighted oriented graphs,} Comm. Algebra  (2024), no.~4, 1672-1685.

\bibitem{gs} D.~R. Grayson , M.~E. Stillman, {\em Macaulay2, a software system for research in algebraic geometry.} Available at \url{http://www.math.uiuc.edu/Macaulay2/}.

\bibitem{hjkn23} H. T. H\`a, A. V. Jayanthan, A. Kumar, H. D. Nguyen, {\em Binomial expansion for saturated and symbolic powers of sums of ideals,} J. Algebra {\bf 620} (2023), 690–710.

\bibitem{hht07}  J. Herzog, T. Hibi and N. V. Trung, {\em Symbolic powers of monomial ideals and vertex cover algebras},  Adv. Math. \textbf{210}  (2007), no.~1, 304-322.

\bibitem{hlmrv19} H. T. H\`a, K. N. Lin, S. Morey, E. Reyes and R. H. Villarreal, {\em Edge ideals of oriented graphs}, Internat. J. Algebra Comput. {\bf 29} (2019), no.~3, 535--559.

\bibitem{htt16} H. T. H\`a, N. V. Trung\ and\ T. N. Trung, {\em Depth and regularity of powers of sums of ideals,} Math. Z. {\bf 282} (2016), no.~3-4, 819--838.

\bibitem{K00} V. Kodiyalam, {\em Asymptotic behaviour of Castelnuovo-Mumford regularity}, Proc. Amer. Math. Soc. {\bf 128} (2000), no.~2, 407--411.

\bibitem{k23} M. Kumar, and R. Nanduri, {\em Regularity of powers of edge ideals of Cohen–Macaulay weighted oriented forests},  J. Algebraic Combin.  {\bf 58}, (2023), no.~3, 867–893. 

\bibitem{mp21} M. Mandal\ and\ D.~K. Pradhan, {\em Symbolic powers in weighted oriented graphs}, Internat. J. Algebra Comput. {\bf 31} (2021), no.~3, 533--549.

\bibitem{mp23} M. Mandal\ and\ D.~K. Pradhan, {\em Properties of symbolic powers of edge ideals of weighted oriented graphs}, Internat. J. Algebra Comput. {\bf 33} (2023), no.~1, 927--951.

\bibitem{m22} N. C. Minh, L. D. Nam , T. D. Phong,
P. T. Thuy, and T. Vu, { \em Comparison between regularity of small symbolic powers and ordinary powers of an edge ideal}, J. Combin. Theory Ser. A {\bf 190} (2022), Paper No. 105621, 30 pp.  

\bibitem{nv19} H. D. Nguyen and T. Vu, {\em Powers of sums and their homological invariants}, J. Pure Appl. Algebra {\bf 223} (2019), no.~7, 3081–3111.

\bibitem{ohh00} H. Ohsugi, J. Herzog and T. Hibi, {\em Combinatorial pure subrings}, Osaka J. Math {\bf 37} (2000), no.~3, 745-757. 

\bibitem{prt19} Y. Pitones, E. Reyes\ and\ J. Toledo, {\em Monomial ideals of weighted oriented graphs}, Electron. J. Combin. {\bf 26} (2019), no.~3, Paper No. 3.44, 18 pp.

\bibitem{s08} S. Sullivant, Combinatorial symbolic powers, J. Algebra {\bf 319} (2008), no.~1, 115--142.

\bibitem{t05} Y. Takayama, Combinatorial characterizations of generalized Cohen-Macaulay monomial ideals, Bull. Math. Soc. Sci. Math. Roumanie (N.S.) {\bf 48(96)} (2005), no.~3, 327--344.

\bibitem{x21} L. Xu, G. Zhu, H. Wang and J. Zhang, {\em Projective dimension and regularity of powers of edge ideals of vertex-weighted rooted forests}, Bull. Malays. Math. Sci. Soc. {\bf 44} (2021), no.~4, 2215--2233.
	 	
\end{thebibliography}
\end{document}